%% file: main.tex
\documentclass[preprint,12pt]{elsarticle}
\input{commands.tex}

\usetikzlibrary{decorations.pathreplacing}
\journal{..}
\begin{document}
	\begin{frontmatter}

\title{Attention-based  hybrid solvers for linear equations that are geometry aware
}

\author{Idan Versano} 
 \ead{idanversano@tauex.tau.ac.il}
\author{Eli Turkel}
 \ead{turkel@tauex.tau.ac.il}
\affiliation{organization={School of Mathematical Sciences, Tel-Aviv University},
            city={ Tel-Aviv},
            postcode={6997801},
            country={Israel}}
	\begin{abstract}
We present a novel architecture for learning geometry-aware preconditioners for linear partial differential equations (PDEs).
We show that a deep operator network (Deeponet) can be trained on a simple geometry and remain a robust preconditioner for problems defined by different geometries without further fine-tuning or additional data mining.
We demonstrate our method for the Helmholtz equation, which is used to solve problems in electromagnetics and acoustics; the Helmholtz equation is not positive definite, and with absorbing boundary conditions, it is not symmetric.
	\end{abstract}


\end{frontmatter}


	\input{introduction.tex}

\input{section2.tex}

\section{Conclusions, Remarks and further research}
We have constructed a novel architecture for preconditioning linear equations, which is aware of the geometry of the problem. We showed evidence by numerical examples
that our network remains a robust preconditioner for problems with different geometries without any fine-tuning procedure and outperforms existing SOTA methods for tackling similar problems. Another advantage of our proposed network  is that is non-convolutional and therefore can be generalized to non-rectangular grids and point-clouds.

We demonstrated our network on a complex Helmholtz equation with complex boundary conditions, which are hard to solve due to non-definiteness or non-symmetry, and on an obstacle problem.

    



In future work, we intend 
to generalize our network to include more attention layers and, with the availability of more advanced hardware, to train the network on a larger number of geometries.

Another possible extension is to modify the neural network so that it is invariant under domain translations, as the PDE is.

\section*{Funding}
The work was supported by the US--Israel Bi-national Science Foundation (BSF) under grant \#~2020128.
\section*{Declaration of competing interest}
The authors declare that they have no known competing financial interests or personal relationships that could have appeared to influence the work reported in this paper.
\section*{Data availability }
\url{
https://github.com/idanv87/convnet.git
}
\input{appendix.tex}

\end{document}

%% file: commands.tex
\usepackage{tabularx}

\usepackage{amsfonts,amssymb,amsmath,amsopn,amsthm}
\usepackage{amsbsy,dsfont}
\usepackage{amsxtra, mathrsfs}
\usepackage{algorithm}
\usepackage{algpseudocode}
\usepackage{graphicx}
\graphicspath{ {./images/} }
\usepackage{float}
\usepackage{colordvi}
\usepackage[usenames,dvipsnames]{color}
\usepackage{nicematrix}
\usepackage{tikz}
	\usepackage{multirow}
\usetikzlibrary{positioning}
\usetikzlibrary{shapes.misc,shadows}
\usepackage{xspace}
\usepackage{lipsum}
\usepackage[many]{tcolorbox}

\usepackage{a4wide,wrapfig,caption,subcaption,epstopdf}

 \numberwithin{equation}{section}

\newcommand{\N}{\mathcal{N}}

\newtcolorbox{leftbrace}{%
	enhanced jigsaw,
	breakable, 
	frame hidden, 
	overlay={%
		\draw [
		decoration={brace,amplitude=0.5em},
		decorate,
		ultra thick,
		]
		(frame.south west)--(frame.north west);
	},
	parbox=false,
}

\usepackage[utf8]{inputenc}
\usepackage{enumitem}

\usepackage{verbatim}

\usepackage{lineno}
\usepackage{textcomp}
\usepackage{mathtools,hyperref}
\usepackage{cleveref}

\usepackage{setspace,esint}
\usepackage{enumitem}

\newcommand{\Hmm}[1]{\leavevmode{\marginpar{\tiny%
			$\hbox to 0mm{\hspace*{-0.5mm}$\leftarrow$\hss}%
			\vcenter{\vrule depth 0.1mm height 0.1mm width \the\marginparwidth}%
			\hbox to
			0mm{\hss$\rightarrow$\hspace*{-0.5mm}}$\\\relax\raggedright #1}}}

\newcommand{\C}{\mathbb{C}}

\newcommand{\R}{\mathbb{R}}

\newtheorem{theorem}{Theorem}[section]

\newtheorem{remark}[theorem]{Remark}
\newtheorem{rem}[theorem]{Remark}

\theoremstyle{definition}

\numberwithin{equation}{section}

\newcommand{\RN}[1]{%
	\textup{\uppercase\expandafter{\romannumeral#1}}%
}

\newcommand{\be}{\begin{equation}}
	\newcommand{\ee}{\end{equation}}
\newcommand{\bea}{\begin{eqnarray}}
	\newcommand{\eea}{\end{eqnarray}}
\newcommand{\bean}{\begin{eqnarray*}}
	\newcommand{\eean}{\end{eqnarray*}}

\newlength{\wex}  \newlength{\hex}
\newcommand{\ass}[1]{Let Assumptions~\ref{assump1} hold  in a bounded Lipschitz domain $\Gw$}

\def\Gw{\Omega}

\usepackage{MnSymbol,bbding,pifont}

\usepackage{multirow}
\usepackage{array}

%% file: introduction.tex
	\section{Introduction}

Solving linear PDE problems frequently involves solving a linear system of equations using iterative algorithms. A central aspect in designing such iterative solvers for a linear system of equations $Au = f$ is finding an approximation of $A^{-1}$, referred to as a preconditioner. The simplest method to utilize a preconditioner $\hat{A} \sim A^{-1}$ is through iterative refinement, also known as Richardson iteration or deferred correction:

\begin{equation}\label{eq:precond}
	u^{i+1} =(1-\theta) u^{i} + \theta \hat{A}(f  -  Au^{i}), \ 0\leq \theta\leq 1. 
\end{equation}
One step of the Richardson iteration method only requires one multiplication of a vector by the matrix $A$, making it suitable for large systems of equations. 

The Helmholtz equation is given by
\begin{equation}\label{eq:helm}
	L u = \Delta u+k^2u=f
\end{equation}
For sufficiently large $k>0$, the operator $L$ is not positive definite  \cite{ernst}. Furthermore, if absorbing boundary conditions are imposed, $L$ is no longer symmetric. It is well known that Richardson iteration with standard solvers like Gauss-Seidel, Jacobi, or geometric multigrid diverge when applied to non-symmetric or non-positive definite operators.
In such cases,  more involved techniques are required, e.g., domain decomposition methods \cite{dolean,gander3, gander4}, Krylov-type methods, \cite{ernst, gander2} including generalized minimal residual method (GMRES) and its variants  \cite{baker, saad,zou}, GMRES with shifted Laplacian preconditioning \cite{gander1}.
However, computational cost and memory storage are more expensive for these methods.

A possible remedy that still allows the use \eqref{eq:richardson} is finding a non-linear preconditioner, $\N$, using a neural network architecture named Deeponet \cite{lu}, which has been experimentally proven to be a robust architecture for learning operators between function spaces, and in particular, for learning solutions of partial differential equations. 
The classical  Deeponet receives as an input $y\in \Gw$ and a function $f$ defined on $\Gw$ and returns $((\Delta+k^2)^{-1}f)(y)$.
The representation of $f$ by finite-dimensional $x\in \R^n$ can be done most straightforwardly by evaluating $f$ on $n$ points in $\Gw$.
The classical Deeponet admits the following general form:
Let $f_m:=(f(x_1),..f(x_m))$ be the values of $f$ over fixed $m$ points in $\Gw$.
Let $\mathcal{B}$ and $\mathcal{T}$ be neural networks such that 
$$
\mathcal{T}:\R^n\to \R^p, \qquad \mathcal{B}:\R^m\to \R^p.
$$ 
Then
\begin{equation}\label{eq:deeponet}
	\mathcal{D}(y,f)=\sum_{j=1}^p\mathcal{T}(y)_j\cdot \mathcal{B}(f_m)_j
\end{equation}
The network $\mathcal{T}$ is named the {\bf trunk} and $\mathcal{B}$ is named the {\bf branch}.
The network $\N$ is trained over data samples such that the loss function   
$$\frac{\|N(y,f)-u(y)\|_{L^2}^2}{|u(y)|^2}$$ is minimized over all samples.
The structure of \eqref{eq:deeponet}  is a classical structure of data-driven recommender systems \cite{zhang1}.

One of the drawbacks of Deeponet is that the numerical accuracy of solutions for a linear system is not fully understood, and reducing the error requires improved architectures \cite{zhu}. In general, solutions with errors of machine precision are not possible. 
To overcome the latter difficulty, a novel hybrid, iterative, numerical, and transferable solver (HINTS) has been recently developed \cite{zhang}. The hybrid solver uses Deeponet as a preconditioner and is integrated inside existing numerical algorithms (e.g., Jacobi, Gauss-Seidel, multigrid) to obtain higher convergence rates for solutions of linear systems arising from discretizations of PDEs up to machine precision.
Classical preconditioners reduce the residual's high-frequency mode.
HINTS balances the convergence behavior across the spectrum of eigenmodes by utilizing the spectral bias of Deeponet, resulting in a uniform convergence rate, giving rise to an efficient iterative solver. This assertion has been recently further investigated in \cite{kopanicakova}, and the robustness of Deeponet-based preconditioners for iterative schemes of linear systems has been shown. We emphasize that Deeponet is trained over a low-resolution grid but provides robust preconditioners for higher resolutions without further training.

A drawback of Deeponets, and as a result also of HINTS, is that it is not geometry aware. Hence,  Deeponet trained on a reference domain might not be efficient for problems on different domains, which are not only small perturbations of the reference domain \cite{kahana}.
The direct solution is modifying Deeponet to train using data from different geometries and multiple domains \cite{he, he2, he3, kashefi, yin, gladstone}.
Yet, as data generation for thousands of domains is time-consuming, we aim to avoid this approach, which might be overkill when using Deeponet as a preconditioner in HINTS.

In \cite{kahana}, a state-of-the-art (SOTA) vanilla Deeponet  has been demonstrated to
operate robustly when used inside HINTS on different domains by extending $f$ to be $0$ for points outside the required domain. 
Though this approach does not require additional training on multiple domains,
it fails to converge in multiple examples we will demonstrate, and training on multiple domains is necessary.

Another potential drawback of the vanilla Deeponet  is its convolutional neural networks (CNN).
Though convolutional neural networks have obtained remarkable capability for studying images and functions on rectangular grids,
in the problems studies by \cite{he, he2, he3, kashefi, yin, gladstone}, it becomes impractical for several 2D/3D problems and several networks excluding convolutional operations were suggested.

In this work, we strive to use HINTS,  also obtained from a {\bf non-convolutional and geometry-aware} Deeponet,  to improve existing iterative schemes to solve problems on different domains even though the training is only on simple geometries. 
We avoid using CNN operations for two reasons:
The first, is that we would like to avoid any artificial extension of function outside its domain.
The second, we wish to provide network that could be used on problems where convolutional operations are impractical.


\subsection{Problem formulation and novelties}
.

We consider the following problem
\begin{equation}\label{eq:sub-problem}
	\begin{cases}
		\Delta u+k^2 u=f &  \Gw\\
		\frac{\partial u}{\partial \vec{\nu}}+\sqrt{-1}ku=0 & \partial \Gw\\
	\end{cases}
\end{equation}
where  $\Gw\subset [0,1]^2$ is a piecewise smooth domain. Again, we stress that this problem is nondefinite for sufficiently large $k$ and nonsymmetric. See \cite[Chapter 35]{ern} for well-posedness.
Standard second-order finite-difference discretization of \eqref{eq:sub-problem} leads to a system of linear equations of the form 
\begin{equation}\label{eq:Axy}
	A{\bf X}={\bf Y}.
\end{equation}
We consider the following  Richardson-type iterations to solve \eqref{eq:Axy}.
\begin{equation}\label{eq:richardson}
	{\bf X}^{n+1}={\bf X}^n+P_n^{-1}r^n
\end{equation}
where $r^n=-A{\bf X}^n+{\bf Y}$ and $P_n^{-1}$ is some preconditioner.
The matrix $ A+A^*$ with $A$ given by  \eqref{eq:sub-problem} can be non-definite and, in general, cheap methods such as Jacobi/Gauss-Seidel diverge.
As noted by \cite{zhang}, $P_n^{-1}$ can be obtained using a neural network $\mathcal{D}$, estimating the inverse operator of \eqref{eq:sub-problem} and lead to convergence of \eqref{eq:richardson} by using the hybrid method 
\begin{equation}\label{eq:hints}
	\begin{cases}
		u^{n+1}=u^n +\mathcal{D}(r^n) & n \ mod \ J=0\\
		u^{n+1}=\mathcal{H}u^n& else
	\end{cases}
\end{equation}
$J\in \mathbb{N}$, and $\mathcal{H}$ is any classical iterative update rule  (e.g., Gauss-Seidel, GMRES).
Our goal is to replace $\mathcal{D}$ by   a neural network $\N$  such that the
resulting scheme  \eqref{eq:hints} applies to families of geometries that are not necessarily "similar"  to the reference geometry used in training without further fine-tuning of the neural network.

\begin{rem}\label{rem:rem1}
	{ \em 
	Deeponet receives as an input a function, and therefore can handle functions defined on different grids, when  the input function is sampled over a finite number of points. 
	By interpolating the input function, Deeponet trained on a certain grid can still handle input data over a coarser or finer grid.  When using coarser/finer grids that contain grid points used for training, the interpolation process is immediate as it requires only choosing a subset of the sample points. Therefore, when one uses Deeponet on different grids, it is less expensive to include grids, which include the grid used for training.
}
\end{rem}
Generating training data requires solving \eqref{eq:sub-problem}.
So, generating data for multiple geometries can be very expensive and challenging when dealing with complex geometries.
Our goal is to generate data on simple geometries only (e.g. rectangles), and then test \eqref{eq:my_hints} on different geometries without any further fine-tuning for different geometry (cf. \cite{goswami}).
As the Neural network takes inputs of fixed sizes, one has to account for dealing with domains with different amounts of grid points. Moreover, one has to deal with how to order the points inside each domain. 
Another problem is how one can evaluate $\N$ on $\Gw$ by computations that "ignore" all points outside $\Gw$, since the non-homogeneous term $f$ in \eqref{eq:sub-problem} is not defined outside $\Gw$. 

\subsection{Main novelties}
We present a new architecture for learning the solutions of \eqref{eq:sub-problem}, which accounts for the domain's geometry.
Our model is a novel improvement to the vanilla Deeponet as it significantly improves the geometry transferability, in case vanilla model might fail. To be more specific, our network is trained on data generated from solutions to  \eqref{eq:sub-problem} on simple reference geometries and still produces robust iterative schemes for different sub-geometries using \eqref{eq:my_hints}, without additional training or fine-tuning (cf. \cite{goswami}).
We examine scenarios where Hints with vanilla deeponet diverge while our network converges.

While common layers in Deeponet are mainly linear and convolutional layers, we also use masked self-attention layers \cite{vaswani}.
The masked self-attention layers enable us, during evaluation on arbitrary geometry $\Gw$ with an arbitrary number of grid points, to include only the points that belong to $\Gw$. This approach does not require an artificial extension of $f$ outside $\Gw$ to fit the input size of the network(cf.  \cite{kahana}).
We observed that when narrow cracks  appear in the tested geometries, the masking effect become significantly noticed and outperforms vanilla model. 

We compare our approaches with the vanilla Deeponet, and naive extension of the vanilla Deeponet to a geometry-aware Deeponet 
with CNN layers.
The latter model is presented to demonstrate that 
our problem is not only a trivial extension of the vanilla case.

To our best knowledge, this is the first attempt to learn  preconditioners using geometry aware neural-networks, and integrating them  with classical iterative solvers.

\subsection{Notations and symbols}\label{subsec:Notations}
We use the following symbols and notations.
\begin{itemize}
	\item For $v\in \R^n$, std(v) is the statistical  standard deviation.
	\item For $v\in \C^n$, std(v)=$std(Real(v)$+$\sqrt{-1} std(Imag(v))$.
	\item For $v,u\in \R^n$, $v\odot u\in \R^n$ is the pointwise multiplication product of $v,u$. 
 \item $h=1/N$ denotes the grid size in a spatial two-dimensional discretization of the unit square $[0,1]^2$ for $N^2$ points.
\item $J$ is the skip factor in \eqref{eq:my_hints}
\item $k$ is a positive constant in the Helmholtz 
equation \eqref{eq:helm}.
\item GMRES corresponds to GMRES algorithm (Algorithm 3 in \cite{saad}).
\item Hints-GS(J) is application of  \eqref{eq:my_hints} where $\mathcal{H}$ is a single application of Gauss-Seidel updating rule.
\item Hints-GMRES(m,J) is application of  \eqref{eq:my_hints} where $\mathcal{H}$ corresponds to  $m$ iterations of GMRES algorithm. 

\end{itemize}

%% file: section2.tex
\newcommand{\setVariableN}{\def\n{80}}

\setVariableN

\section{Geometry aware Deeponet}
In this section we describe in details our proposed model.
\subsection{The model}
We generate a uniform grid on $[0,1]^2$ with $15 \times 15$ points. We flatten the grid points to obtain
a vector  $\vec{d}\in \R^{225\times 2}$.

{\bf Model's inputs:}
\begin{enumerate}
\item $\vec{d}$.
\item 
A domain $\Gw\subset [0,1]^2 $ represented by vector $O^{\Gw}\in \R^{225}$, where,
$$
O_i^{\Gw}=\text{dist}((d_{i0},d_{i1}), \partial \Gw),
$$
and $\text{dist}$ is the signed distance function.
\item A real function $f:\Gw \to \R$, represented by the complex vector $f^{\Gw}\in \C^{225}$, where
$$
f^{\Gw}_i=
\begin{cases}
    f(d_{i0},d_{i1}),&  (d_{i0},d_{i1})\in \Gw \\
    0, & otherwise.
\end{cases}
$$
\item A masking matrix $M^{\Gw}\in \R^{225\times 225}$  defined as follows.
$$
M^{\Gw}=
\begin{cases}
	0 & (d_{i0},d_{i1})\in \Gw_n \\
	-\infty & otherwise.
\end{cases}
$$
\item $y\in \Gw$.
\end{enumerate}

{\bf Model's output}: 
a complex number reresenting  the solution to \eqref{eq:helm}  evaluated at $y$.

\subsubsection{Model's layers}
First, we define the M-attention operation of $v\in \R^N$ and $M^{\Gw}$ to be 
$$
\text{M-attention}(v,M^{\Gw}):=
\text{softmax}\left (\frac{vv^T+M^{\Gw}}{\sqrt{N}}\right )v
$$
\cite{vaswani}.
The function  $\text{softmax}:\R^k\to (0,1)^k$ is defined by 
$$
\text{softmax}({\bf Z})_i =\frac{\exp({\bf Z_i})}{\sum_{i=1}^k \exp({\bf Z_i})}.
$$
Notice that the M-attention does not depend on the values of $f$ outside $\Gw$.
If the masking matrix $M^{\Gw}$ is replaced with a zero matrix, then M-attention accounts for a trivial extension of $f$ outside $\Gw$ and, therefore, depends on the extension one chooses for $f$ outside $\Gw$.

Since the output of the model should be a complex number we first define two Deeponets,
 $\mathcal{N}_{real}$,
$\mathcal{N}_{imag}$  having the  structure \eqref{eq:deeponet}  with the following branch and trunk:

Branch:
$$
\mathcal{B}:\R^{225\times 4} \to \R^{80}
$$
defined by the operation of
 $M^{\Gw}$-attention on the matrix  generated by concatenating $f^{\Gw},\vec{d}, O^{\Gw}$, and followed by a dense NN [225,200,100,80] with a Tanh activation.
 
Trunk:
$$
\mathcal{T}:\R^{3} \to \R^{80}
$$
defined by the operation of a  dense neural network (NN) with layers [3,200,100,80] and $\sin$ activation function (\cite{he}), on  the vector  generated by concatenating $y$ and 
$\text{dist}(y,\partial \Gw)$.

We construct a neural Network 
\begin{equation}\label{eq:Nmasked}
\N^{masked}=\N_{real}+\sqrt{-1}\N_{imag}.
\end{equation}
We also define  a second  network 
\begin{equation}\label{eq:Nnonmasked}
\N^{non-masked}=\N_{real}+\sqrt{-1}\N_{imag},
\end{equation}
where the masking matrix $M^{\Gw}$ is replaced by zero matrix, and an input function $f$ is extended by zeros outside its domain.

We  define an updating rule for the hybrid method called
HINTS,
as follows:
\begin{equation}\label{eq:my_hints}
	\begin{cases}
		u^{n+1}=u^n +\frac{std(r^n)}{\alpha}\cdot \N\left ( \alpha \frac{r^n}{std(r^n)}\right ) 
		& n \ mod \ J=0\\
		u^{n+1}=\mathcal{H}u^n& otherwise
	\end{cases}
\end{equation} 
where $J\in \mathbb{N}$  and $\alpha>0 $ are (non-trained) hyper-parameters to be selected by the user, and $n\in \mathbb{N}$ is the iteration number.
The normalization factor $\frac{\alpha}{std(r^n)}$ ensures the input of the network $\N$ will belong to the distribution of the training data.

\begin{remark}
{\em 
We have defined our networks to consider as inputs only real functions.
In order to apply $\mathcal{N}^{masked}$ ($\mathcal{N}^{non-masked}$) on complex functions one can apply the same network separately on the real and imaginary part of $f$.
}
\end{remark}
\subsection{Model's training data}\label{subsec:training}
\begin{figure}[!t]
	\centering
	\subfloat[\centering ]{{\includegraphics[width=12cm, scale=1]{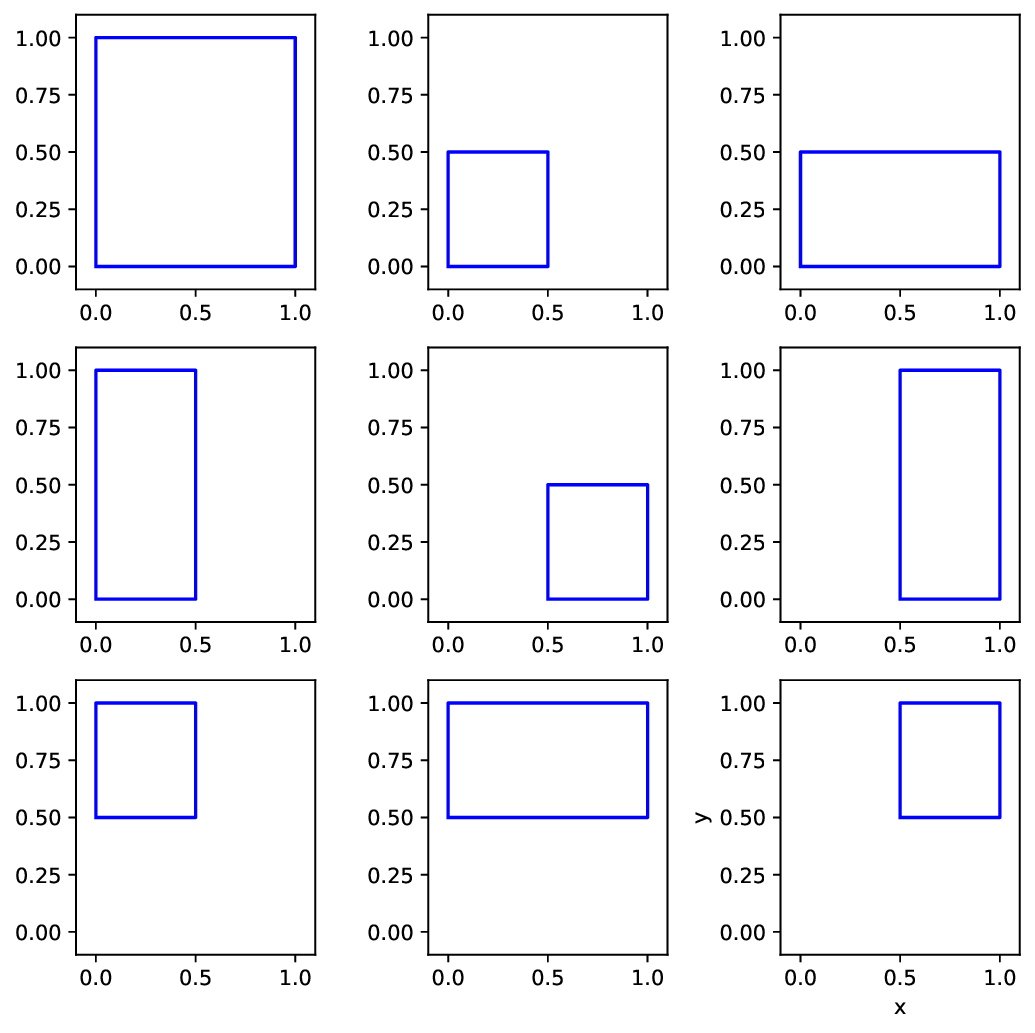} }}
	\caption{Geometries used for training }.
	\label{fig:trained_geometries}
\end{figure}

We fix $k^2=21$.
Then, we discretize each of  the training geometries  in Figure  \ref{fig:trained_geometries}   with a uniform grid of  grid size $h=1/14$.
The training data for $\N$ is obtained by solving \eqref{eq:sub-problem} on the geometries in Figure \ref{fig:trained_geometries} using standard second-order finite difference method,    with different 300 functions $f$, taken to be  Gaussian random fields over $\R$ with zero mean and covariance kernel 
\begin{equation}\label{eq:cov}
cov(f({\bf x}_i), f({\bf x}_j))=\sigma\exp \left (- \frac{1}{2}\frac{\| {\bf x}_i- {\bf x}_j\|^2}{l^2} \right ),
\end{equation}
where 
$$
l=\sigma=0.1.
$$
\subsection{Network's Complexity}
Every $J$ iterations of the classical iterative method is followed by one application of $\mathcal{N}$. 
As the network's input is a vector of fixed dimension, each evaluation of the Network on a single data point is of constant computational and memory complexity.
As a result, for a sparse (and band-limited) linear system of order $N$, the complexity of each iteration involving $\mathcal{N}$ in \eqref{eq:my_hints} is $O(N)$ for both memory and computational complexity.

\newpage
\section{Numerical simulations}\label{sec:simulations}

We compare the performance of HINTS \eqref{eq:my_hints} for four different models:
\begin{enumerate}
    \item $\mathcal{N}^{masked}$ is given by            \eqref{eq:Nmasked} and was trained on geometries in Figure \ref{fig:trained_geometries}.
    \item $\mathcal{N}^{non-masked}$ is given by \eqref{eq:Nnonmasked} and was trained on geometries in Figure \ref{fig:trained_geometries}.
    \item vanilla Deeponet taken from \cite{zhang} (see  \ref{sec:app}) and was trained on $[0,1]^2$.
    \item ga-vanilla is a naive implementation of geometry-aware Deeponet using CNN (see  \ref{sec:app}) that was trained on the geometries in Figure \ref{fig:trained_geometries}.
        
\end{enumerate}

\begin{figure}[!t]
	\centering
	\subfloat[\centering ]{{\includegraphics[width=15cm, scale=1]{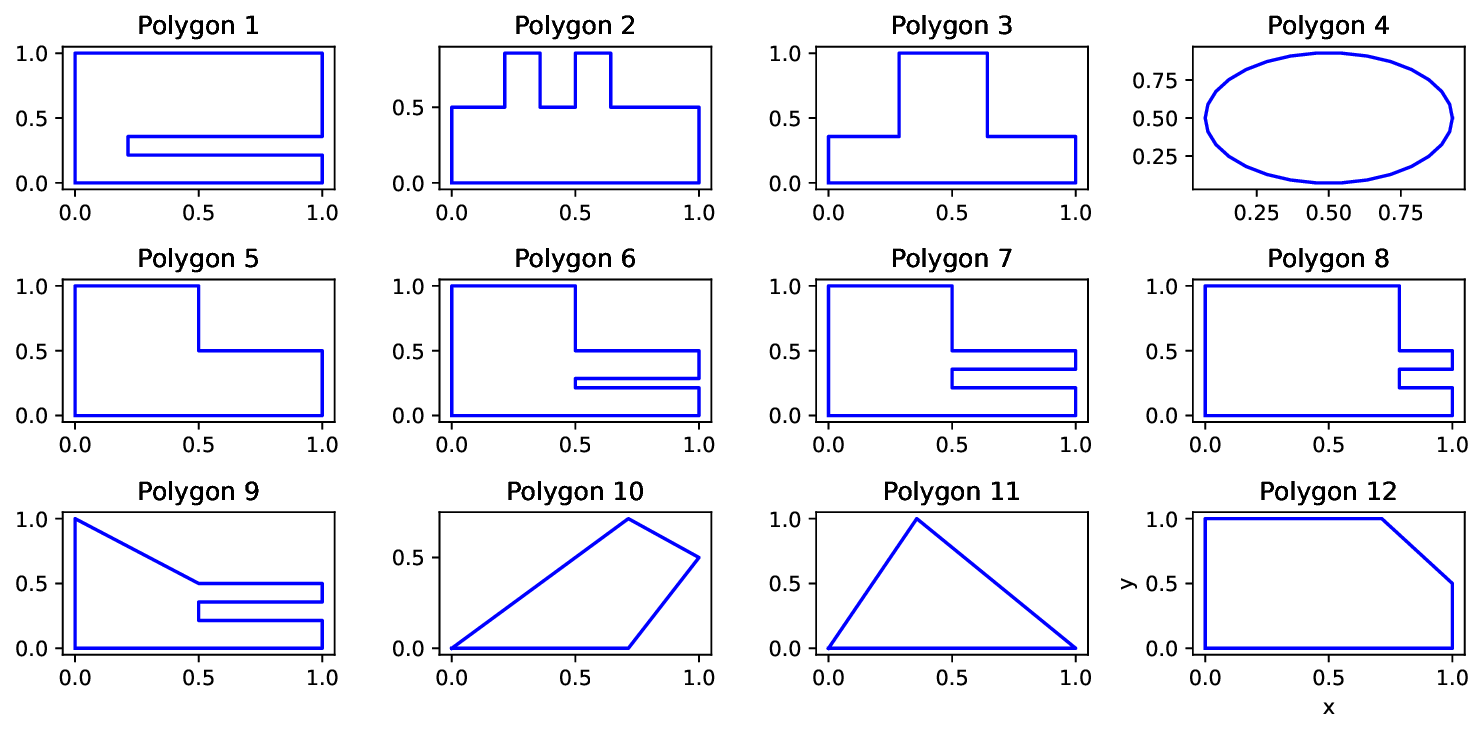} }}
	\caption{Tested domains $\Gw$ used in Tables \ref{tab:exp}, \ref{tab:exp2} and \ref{tab:expgmews} for solving \eqref{eq:helm}.
 }.
	\label{fig:tested_geometries}
\end{figure}

\begin{table}[htbp]
	\centering
	\begin{tabular}{|c|c|c|c|}
		\hline
		$h$ & Method& Model&iter. number\\ 
		\hline
\multicolumn{4}{|c|}{geometry 1} \\
		\hline
1/57 & H-GS(J=40) &$\mathcal{N}^{masked}$&4877  \\
1/57 & H-GS(J=40) &ga-vanilla&diverge  \\
1/57 & H-GS(J=40) &vanilla&diverge  \\
1/57 & H-GS(J=40) &$\mathcal{N}^{non-masked}$&diverge  \\

\hline
\multicolumn{4}{|c|}{geometry 2} \\
		\hline
1/57 & H-GS(J=60) &$\mathcal{N}^{masked}$&12000  \\
1/57 & H-GS(J=60) &ga-vanilla&diverge  \\
1/57 & H-GS(J=60) &vanilla&diverge  \\
1/57 & H-GS(J=60) &$\mathcal{N}^{non-masked}$& 16000 \\

\hline
\multicolumn{4}{|c|}{geometry 3} \\
		\hline
1/57 & H-GS(J=40) &$\mathcal{N}^{masked}$&$2157$  \\
1/57 & H-GS(J=40) &ga-vanilla&1616  \\
1/57 & H-GS(J=40) &vanilla&$1325$  \\
1/57 & H-GS(J=40) &$\mathcal{N}^{non-masked}$&$1381$  \\

\hline
\multicolumn{4}{|c|}{geometry 4} \\
		\hline
1/57 & H-GS(J=40) &$\mathcal{N}^{masked}$&$1519$ \\
1/57 & H-GS(J=40) &ga-vanilla&diverge  \\
1/57 & H-GS(J=40) &vanilla&$813$  \\
1/57 & H-GS(J=40) &$\mathcal{N}^{non-masked}$&$764$  \\

\hline
\multicolumn{4}{|c|}{geometry 5} \\
		\hline
1/57 & H-GS(J=40) &$\mathcal{N}^{masked}$& $1620\pm 170$ \\
1/57 & H-GS(J=40) &ga-vanilla&diverge  \\
1/57 & H-GS(J=40) &vanilla&$1064$  \\
1/57 & H-GS(J=40) &$\mathcal{N}^{non-masked}$&$842$  \\

\hline
\multicolumn{4}{|c|}{geometry 6} \\
		\hline
1/57 & H-GS(J=40) &$\mathcal{N}^{masked}$&6300 \\
1/57 & H-GS(J=40) &ga-vanilla&diverge  \\
1/57 & H-GS(J=40) &vanilla&$>$20000  \\
1/57 & H-GS(J=40) &$\mathcal{N}^{non-masked}$&7200  \\

\hline
\multicolumn{4}{|c|}{geometry 7} \\
		\hline
1/57 & H-GS(J=40) &$\mathcal{N}^{masked}$&10800  \\
1/57 & H-GS(J=40) &ga-vanilla&diverge  \\
1/57 & H-GS(J=40) &vanilla&diverge  \\
1/57 & H-GS(J=40) &$\mathcal{N}^{non-masked}$&8034  \\

\hline
\multicolumn{4}{|c|}{geometry 8} \\
		\hline
1/57 & H-GS(J=40) &$\mathcal{N}^{masked}$&2028  \\
1/57 & H-GS(J=40) &ga-vanilla&2400  \\
1/57 & H-GS(J=40) &vanilla&$1370$  \\
1/57 & H-GS(J=40) &$\mathcal{N}^{non-masked}$&$1617$  \\

\hline

	\end{tabular}
	\caption{Averaged iteration number statistics for the method \eqref{eq:my_hints} with HINTS Gauss-Seidel (H-GS) with the geometries in Figure \ref{fig:tested_geometries}, $k^2=21$, $\alpha=0.3$.
 The classical Gauss-Seidel iterations are not presented as they diverge.}
	\label{tab:exp}
\end{table}

\input{tables.tex}

Our tested geometries for problem \eqref{eq:helm} are presented in Figure \ref{fig:tested_geometries}.
We discretize each tested domain with a higher resolution than the resolution used for training, $N=15$.
We solve \eqref{eq:sub-problem} using \ref{eq:my_hints} with a random function $f$,
sampled from another distribution than the one used for training ( $f$ has been generated as a vector from a normal distribution with mean ten and variance 10.).
Since the neural network operates on the error terms only, the distribution of $f$ does not affect the hybrid scheme's robustness.
We compare our method with both Gauss-Seidel iterations
(see Table \ref{tab:exp})
and GMRES iterations (see Table \ref{tab:expgmews}).
The main advantage of Gauss-Seidel (GS) is that, for sparse systems, it has linear complexity in both memory and computational time. However, GS does not always converge for the non-definite cases in our examples.
We will show that Hints-GS converges even when plain Gauss-Seidel iterations diverge.
GMRES's main advantage is proven convergence for any non-singular matrix, while its major drawback is memory complexity. 
For (not necessarily sparse) linear equations with $N$ equations (i.e., corresponding matrix $A\in \C^{N\times N}$), the iteration complexity of GMRES after $n$ iterations is $O(Nn^2)$ and the memory complexity is $O(N)$.
A possible remedy is the restarted GMRES method, GMRES(m), whose memory complexity reduces to $O(mn)$ independently of the number of iterations but with a possible degradation of convergence rates.
We will show that combining GMRES inside \eqref{eq:my_hints} yields a more efficient method in terms of both time and computational complexity.

It is worth noting that, even when Gauss-Seidel iterations diverge, they are still known to reduce high-frequency errors throughout the iterations. Therefore, Hints-GS highlights the capability of HINTS to effectively control low-frequency modes.

 
 We use the following schemes for comparisons as $\mathcal{H}$ in \eqref{eq:my_hints}:
 GMRES, Gauss-Seidel (GS) , Hints-GS, Hints-GMRES, see subsection \ref{subsec:Notations}.
 Multiple computations were performed, and we show the mean number of iterations required. $k>0$ is the wavenumber in the Helmholtz equation, $h=dx=dy$ is the discretization size, and $J$ is the skip factor in \eqref{eq:my_hints} for how often the neural network $\mathcal{N}$ (Deeponet) is called.
\begin{remark}
{ \em 
One of the highlights of \cite{kahana} was the generalization capability of the vanilla Deeponet on a difficult obstacle problem.
    We also performed experiments  of problem \eqref{eq:helm} on an obstacle problems with $\Gw=[0,1]^2 \setminus [0.25,0.5]^2$ and
    $\Gw=[0,1]^2 \setminus [0.25,0.75]^2$, and observed the superiority of $\mathcal{N}^{masked}$ on the vanilla Deeponet (see Table \ref{tab:expobs}).
    We have also performed additional experiments on exterior domains with different shapes of obstacles  
    (see Table \ref{tab:expmoreobs}).
    }
    \end{remark}
    
    \begin{table}[htbp]
	\centering
	\begin{tabular}{|c|c|c|c|c|}
		\hline
		$h$ & Method& Model&iter. number\\ 
		\hline
  \multicolumn{4}{|c|}{obstacle 1} \\
		\hline
1/57 & Gauss-Seidel &-&diverge  \\  
1/57 & H-GS(J=40) &vanilla&18000 \\
1/57 & H-GS(J=40) &$\mathcal{N}^{masked}$&8000 \\
\hline
  \multicolumn{4}{|c|}{obstacle 2} \\
		\hline
1/57 & Gauss-Seidel &-&10600 \\  
1/57 & H-GS(J=100) &vanilla&$>$20000 \\
1/57 & H-GS(J=100) &$\mathcal{N}^{masked}$&3425 \\
\hline
  \multicolumn{4}{|c|}{obstacle 1} \\
		\hline
1/224 & GMRES &-&21400 \\  
1/224 & H-GMRES(J=15) &vanilla&diverge \\
1/224 & H-GMRES(J=15) &$\mathcal{N}^{masked}$&14100 \\
\hline
  \multicolumn{4}{|c|}{obstacle 2} \\
		\hline
1/224 & GMRES &-&6460 \\  
1/224 & H-GMRES(J=20) &vanilla&12000 \\
1/224 & H-GMRES(J=20) &$\mathcal{N}^{masked}$&8100 \\
\hline
	\end{tabular}
	\caption{Averaged iteration number statistics for the method \eqref{eq:my_hints} with HINTS Gauss-Seidel (H-GS) with the geometries in Figure \ref{fig:fig_obstacle}, $k^2=21$, $\alpha=0.3$.}
	\label{tab:expobs}
\end{table}
\begin{figure}[!t]
	\centering
	\subfloat[\centering ]{{\includegraphics[width=15cm, scale=1]{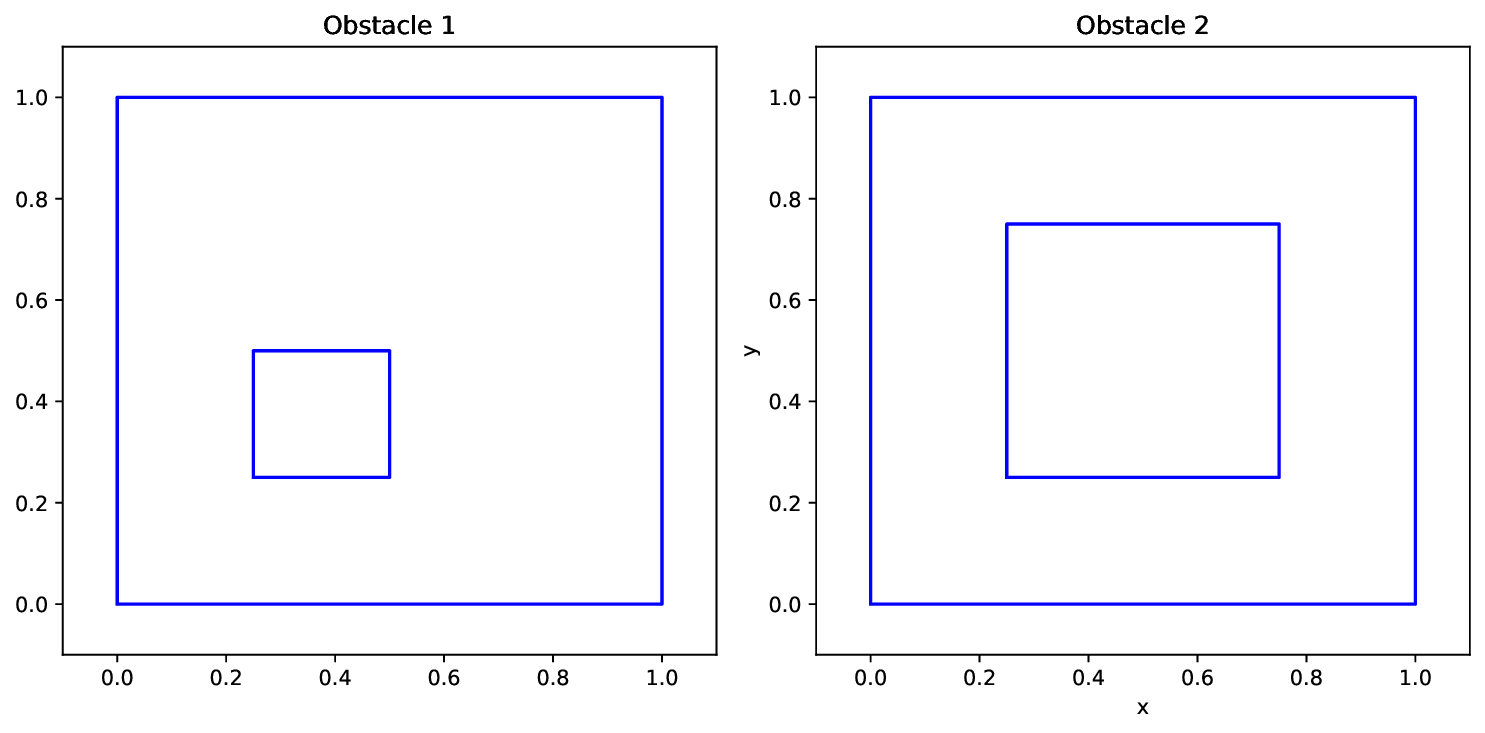} }}
	\caption{Tested exterior domains used in Table \ref{tab:expobs} for solving $\Delta u+k^2 u=f$. The boundary conditions on the inner boundary are u=0.
 The boundary conditions on the outer boundary are $\frac{\partial u}{\partial \vec{\nu}}+\sqrt{-1}k u=0$.
 }
	\label{fig:fig_obstacle}
\end{figure}

    \begin{table}[htbp]
	\centering
	\begin{tabular}{|c|c|c|c|c|}
		\hline
		$h$ & Method& Model&iter. number\\ 
		\hline
  \multicolumn{4}{|c|}{1 obstacle } \\
		\hline
1/57 & Gauss-Seidel &-& diverge \\  
1/57 & H-GS(J=100) &vanilla&16700 \\
1/57 & H-GS(J=100) &$\mathcal{N}^{masked}$&8500 \\
\hline
  \multicolumn{4}{|c|}{2 obstacles } \\
		\hline
1/57 & Gauss-Seidel &-& diverge \\  
1/57 & H-GS(J=100) &vanilla&$>$40000 \\
1/57 & H-GS(J=100) &$\mathcal{N}^{masked}$&14500 \\
\hline
  \multicolumn{4}{|c|}{3 obstacles } \\
		\hline
1/57 & Gauss-Seidel &-& 12500 \\  
1/57 & H-GS(J=200) &vanilla& 8300\\
1/57 & H-GS(J=200) &$\mathcal{N}^{masked}$&8500 \\
\hline
  \multicolumn{4}{|c|}{1 obstacle } \\
		\hline
1/224 & GMRES &-& 32100 \\  
1/224 & H-GMRES(m=30,J=15) &vanilla&10575 \\
1/224 & H-GMRES(m=30,J=15) &$\mathcal{N}^{masked}$& 16480 \\
\hline
  \multicolumn{4}{|c|}{2 obstacles } \\
		\hline
1/224 & GMRES &-&18000  \\  
1/224 & H-GMRES(m=30,J=15) &vanilla&14193 \\
1/224 & H-GMRES(m=30,J=15) &$\mathcal{N}^{masked}$& 15810 \\
\hline
  \multicolumn{4}{|c|}{3 obstacles } \\
		\hline
1/224 & GMRES &-& 8600 \\  
1/224 & H-GMRES(m=30,J=15) &vanilla&13343 \\
1/224 & H-GMRES(m=30,J=15) &$\mathcal{N}^{masked}$& 16717 \\
\hline
	\end{tabular}
	\caption{Averaged iteration number statistics for the method \eqref{eq:my_hints} with HINTS Gauss-Seidel (H-GS) with the geometries in Figure \ref{fig:fig_more_obstacle}, $k^2=21$, $\alpha=0.3$.}
	\label{tab:expmoreobs}
\end{table}
\begin{figure}[!t]
	\centering
	\subfloat[\centering ]{{\includegraphics[width=15cm, scale=1]{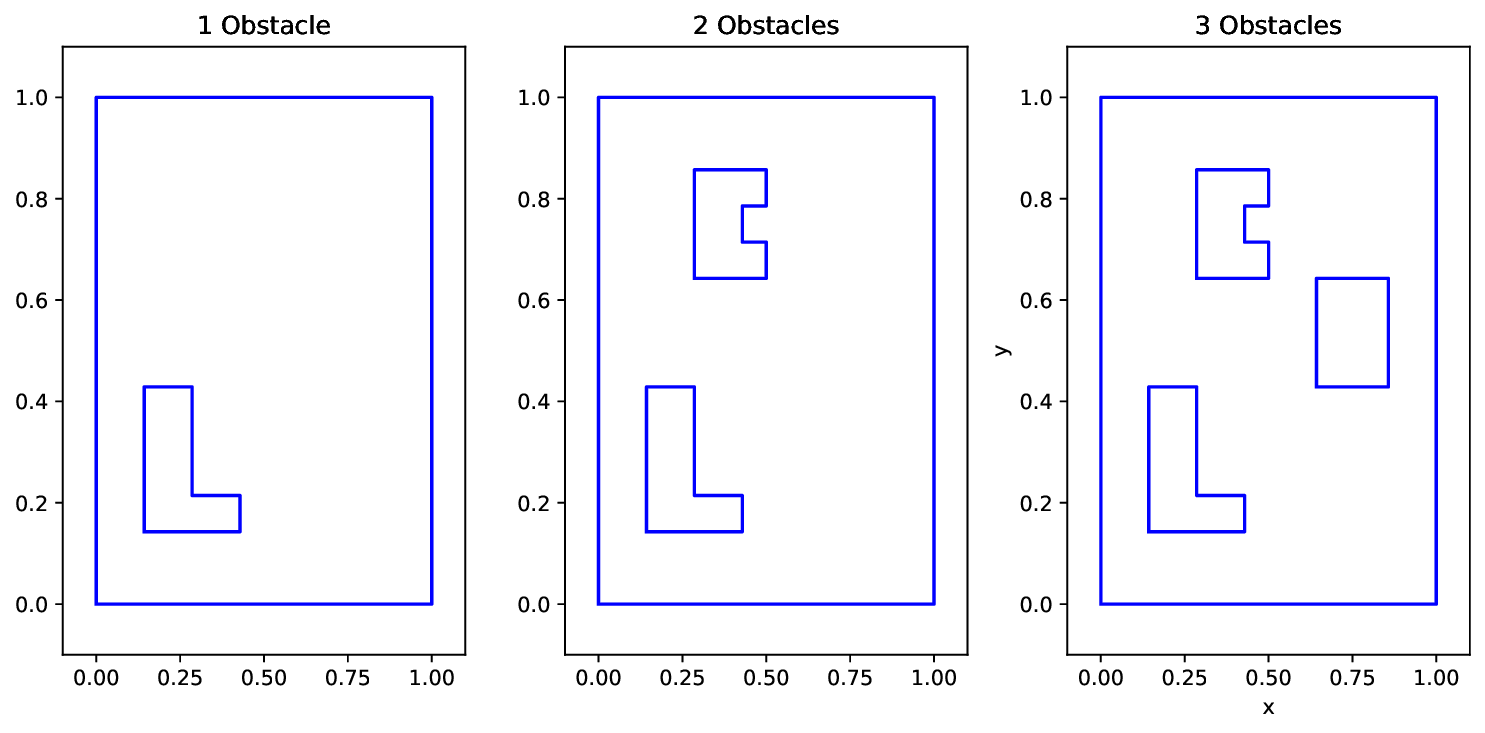} }}
	\caption{Tested exterior domains used in Table \ref{tab:expmoreobs} solving $\Delta u+k^2 u=f$. The boundary conditions on the inner boundaries are u=0.
 The boundary conditions on the outer boundary are $\frac{\partial u}{\partial \vec{\nu}}+\sqrt{-1}k u=0$.
 }
	\label{fig:fig_more_obstacle}
\end{figure}
\newpage

\subsection{Observations from Tables \ref{tab:exp}, \ref{tab:exp2}, \ref{tab:expgmews}, \ref{tab:expobs},\ref{tab:expmoreobs}}, \label{subsec:observations}

We have demonstrated multiple examples with complex geometries for which HINTS with vanilla method fails, while our network succeed. This confirms the necessity of the new models $\mathcal{N}^{masked}, \mathcal{N}^{non-masked}$   we presented.
In particular, in terms of convergence capability, our new models outperforms the vanilla model.
Yet, we mention that the geometry transferability
of the vanilla still perform robustly on many domains including
non-rectangular ones (geometries 3,4,5,8), and in many cases
provide similar or slightly better convergence rates to that of
our the geometry aware models.

In  almost all tests ga-vanilla obtained the worse
results, implying that a naive extension  of the vanilla model does not work and  more sophisticated approach was required.

When the domain's boundary contains abrupt changes (cracks/bumps/hole), geometries 1,2,6,7,9, $\mathcal{N}^{masked}$ performed better then all methods.
In cases for which vanilla outperformed  $\mathcal{N}^{masked}$,
$\mathcal{N}^{non-masked}$  outperformed $\mathcal{N}^{masked}$, with an exception for H-GMRES for geometry 1.

In almost all cases we tested there is correlation between HINTS' performance with different methods (GS,GMRES). Yet, the method used inside HINTS and the corresponding parameters can violate this correlation (e.g.  Table \ref{tab:expmoreobs}). In more details, the experiments with the geometries in Figure \ref{fig:fig_more_obstacle} shows that with GS, $\mathcal{N}^{masked}$ outperforms the vanilla method significantly, whereas for the same problem with GMRES (J=15)
vanilla method is better. Yet, for different parameter $J=5$ H-GMRES with vanilla diverged while it did converged with $\mathcal{N}^{masked}$.

We can tell that the idea of extending  a function outside its domain to be 0, as suggested in \cite{kahana} can indeed provide robust results for different domains as long as the domain does not contain abrupt changes (narrow cracks, holes).
In these problematic scenarios our masked model can overcome this limit.

Based on our results we suggest to a potential user which cannot use CNN operations to hold both 
$\mathcal{N}^{masked}$ and $\mathcal{N}^{non-masked}$.
Further study on unifying these networks is left for a future research.

\begin{remark}
{\em 
    The running times of \eqref{eq:my_hints} can differ for different types of computers. There is room for additional improvement of \eqref{eq:my_hints} as we did not use a GPU and did not exploit a parallel computation of batched data during inference. In our machine H-GS with higher resolutions than $1/57$ is either very slow or diverges and therefore we switched to GMRES to obtain faster CPU time.
    }
\end{remark}

%% file: tables.tex
\begin{table}[htbp]
	\centering
	\begin{tabular}{|c|c|c|c|}
		\hline
		$h$ & Method& Model&iter. number\\ 
		\hline
\multicolumn{4}{|c|}{geometry 9} \\
		\hline
1/57 & H-GS(J=80) &$\mathcal{N}^{masked}$&9700  \\
1/57 & H-GS(J=80) &ga-vanilla&diverge  \\
1/57 & H-GS(J=80) &vanilla&diverge  \\
1/57 & H-GS(J=80) &$\mathcal{N}^{non-masked}$&18500  \\

\hline
\multicolumn{4}{|c|}{geometry 10} \\
		\hline
1/57 & H-GS(J=40) &$\mathcal{N}^{masked}$&2500  \\
1/57 & H-GS(J=40) &ga-vanilla&diverge  \\
1/57 & H-GS(J=40) &vanilla&2450  \\
1/57 & H-GS(J=40) &$\mathcal{N}^{non-masked}$&1055  \\

\hline
\multicolumn{4}{|c|}{geometry 11} \\
		\hline
1/57 & H-GS(J=80) &$\mathcal{N}^{masked}$&2100  \\
1/57 & H-GS(J=80) &ga-vanilla&1770  \\
1/57 & H-GS(J=80) &vanilla&$1055$  \\
1/57 & H-GS(J=80) &$\mathcal{N}^{non-masked}$&1040 \\

\hline
\multicolumn{4}{|c|}{geometry 12} \\
		\hline
1/57 & H-GS(J=40) &$\mathcal{N}^{masked}$&1200  \\
1/57 & H-GS(J=40) &ga-vanilla&1000  \\
1/57 & H-GS(J=40) &vanilla&940  \\
1/57 & H-GS(J=40) &$\mathcal{N}^{non-masked}$&560  \\

\hline
	\end{tabular}
	\caption{Averaged iteration number statistics for the method \eqref{eq:my_hints} with HINTS Gauss-Seidel (H-GS) with the geometries in Figure \ref{fig:tested_geometries}, $k^2=21$, $\alpha=0.3$.
 The classical Gauss-Seidel iterations are not presented as they diverge.}
	\label{tab:exp2}
\end{table}

\begin{table}[htbp]
	\centering
	\begin{tabular}{|c|c|c|c|c|}
		\hline
		$h$ & Method& Model&GMRES iter.num.&Deeponet iter. num.\\ 
		\hline
\multicolumn{5}{|c|}{geometry 1} \\
		\hline
1/224 & GMRES&-&27120&- \\  
1/224 & H-GMRES(m=30,J=10)&$\mathcal{N}^{masked}$&11400&130  \\ 
1/224 &H-GMRES(m=30,J=10)&vanilla&11000&120 \\ 
1/224 &H-GMRES(m=30,J=10)&$\mathcal{N}^{non-masked}$&24797&273 \\ 
\hline

\multicolumn{5}{|c|}{geometry 2} \\
		\hline
1/224 & GMRES&-&19416&-  \\  
1/224 & H-GMRES(m=30,J=10)&$\mathcal{N}^{masked}$&8500& 90 \\ 
1/224 &H-GMRES(m=30,J=10)&vanilla&25800&285 \\ 
1/224 &H-GMRES(m=30,J=10)&$\mathcal{N}^{non-masked}$&7200 &70 \\ 
\hline

\multicolumn{5}{|c|}{geometry 4} \\
		\hline
1/224 & GMRES&-&21570&-  \\  
1/224 & H-GMRES(m=30,J=10)&$\mathcal{N}^{masked}$&5300&57 \\ 
1/224 &H-GMRES(m=30,J=10)&vanilla&3824&42 \\ 
1/224 &H-GMRES(m=30,J=10)&$\mathcal{N}^{non-masked}$&4490&48 \\ 
\hline
\multicolumn{5}{|c|}{geometry 6} \\
		\hline
1/224 & GMRES&-&17816&-  \\  
1/224 & H-GMRES(m=30,J=10)&$\mathcal{N}^{masked}$&8814&96 \\ 
1/224 &H-GMRES(m=30,J=10)&vanilla&13400&147 \\ 
1/224 &H-GMRES(m=30,J=10)&$\mathcal{N}^{non-masked}$&13079&144 \\ 
\hline

\multicolumn{5}{|c|}{geometry 7} \\
		\hline
1/224 & GMRES&-&18426&-  \\  
1/224 & H-GMRES(m=30,J=10)&$\mathcal{N}^{masked}$&8800&96 \\ 
1/224 &H-GMRES(m=30,J=10)&vanilla&16500&190 \\ 
1/224 &H-GMRES(m=30,J=10)&$\mathcal{N}^{non-masked}$&15300&170 \\ 
\hline
\multicolumn{5}{|c|}{geometry 9} \\
		\hline
1/224 & GMRES&-&14616&-  \\  
1/224 & H-GMRES(m=30,J=10)&$\mathcal{N}^{masked}$&8800&100 \\ 
1/224 &H-GMRES(m=30,J=10)&vanilla&14600&160 \\ 
1/224 &H-GMRES(m=30,J=10)&$\mathcal{N}^{non-masked}$&9900&108 \\ 
\hline
\multicolumn{5}{|c|}{geometry 11} \\
		\hline
1/224 & GMRES&-&16257&-  \\  
1/224 & H-GMRES(m=30,J=10)&$\mathcal{N}^{masked}$&5500&60 \\ 
1/224 &H-GMRES(m=30,J=10)&vanilla&4140&45 \\ 
1/224 &H-GMRES(m=30,J=10)&$\mathcal{N}^{non-masked}$&3830&42 \\ 
\hline
	\end{tabular}
	\caption{Averaged iteration number statistics for the method \eqref{eq:my_hints} with HINTS GMRES (H-GMRES) with the geometries in Figure \ref{fig:tested_geometries}, $k^2=21$, $\alpha=0.3$.}
	\label{tab:expgmews}
\end{table}

%% file: appendix.tex
\appendix
\section{}

\subsection{The Networks}\label{sec:app}
For vanilla deeponet, trained on a single domain $[0,1]^2$, we employ the architecture in \cite[page 37]{zhang}.
The branch net is a combination of a 2D CNN (input dimension $15\times 15$, number of channels [1, 40, 60, 100], kernel size $3\times 3$, stride 2) and dense NN (dimension [100, 80, 80]) with Relu activation function. The dimension of the dense trunk network is [2, 80, 80, 80] with activation of Tanh. 

For the naive geometry-aware deeponet (ga-vanilla), we represent the domain as 2D binary image and pad the values of the function $f$ with zeros in all points outside the domain $f$ is defined. We then use CNN as is the vanilla deeponet whose output is concatenated and fed into dense NN (dimension [200, 80, 80]), resulting in the branch.
The dimension of the trunk network is [2, 80, 80, 80] with the activation of Tanh. 

\subsection{Data generation and training }\label{sec:app1}
We fix $k^2=21$ in \eqref{eq:sub-problem}.
The vanilla network  is trained over 300 solutions of \eqref{eq:sub-problem} in $[0,1]^2$ discretized  uniformly  with $h=\frac{1}{15}$. The solutions are obtained by taking $f$ to be  300 Gaussian random fields \eqref{eq:cov}.  The equations are solved using  standard second-order finite difference method.
The network ga-vanilla is trained on the same data used to train $\mathcal{N}^{masked}$, see subsection \ref{subsec:training}.

We trained the networks using Pytorch for 1000 epochs with a batch size of 64, and the initial learning rate 1e-4 decreased by 0.5 at the final epochs.